\documentstyle{amsppt}

\document

\NoRunningHeads

\topmatter

\title \bf On some examples in Symplectic Topology
\endtitle

\footnote ""{Partially supported by Russian Foundation of Basic Research
 (grant N 01-01-00709)}

\author V. Gorbatsevich
\endauthor

\abstract

Article is devoted to the Examples 2 and 3 of the symplectic 
solvable Lie groups $R$ with some special cohomological 
properties, which have been constructed by Benson and Gordon. 
But they are not succeeded in constructing corresponding compact 
forms for symplectic structures on these Lie groups. Recently 
A.Tralle proved that there is no compact form in the Example 3. 
But his proof is rather complicated and uses some very special 
topological methods. 

We propose much more simpler (and purely algebraic) method to prove 
the main result of the Tralle's paper. Moreover we prove that for 
Example 2 there is no compact form too. But it appears that some 
modification  of the construction of the Example 2 gives some other 
example of a solvable Lie group $R^{\prime}$ with the same cohomological 
properties as $R$, but with a compact form. 

\endabstract



\endtopmatter

Let $(M,\omega)$ be a compact symplectic manifold. This article is devoted to some examples, 
which have been constructed in [1], where some problems about K\"{a}hlerian structures on 
solvmanifolds (i.e. homogeneous spaces $R/\Gamma$ of solvable Lie groups $R$ with discrete stationary 
subgroups $\Gamma$) are studied. There are Examples 2 and 3 of the symplectic solvable Lie groups $R$ 
in [1] with some special properties (closely related to the properties of K\"{a}hlerian Lie groups). 
But the authors of [1] are not succeeded in constructing corresponding compact forms $R /\Gamma$
for this symplectic structures on Lie groups. These examples in [1] have been constructed for the 
purpose of illustration of the conditions of the main result (see Theorem 2 there about a structure 
of compact solvmanifolds $R/\Gamma$ which admit a K\"{a}hler structure) of this paper. 
Recently A.Tralle [2] proved that there is no compact form (i.e. a compact solvmanifold 
$R/\Gamma$ for the Lie group $R$) in the Example 3 from [1]. But the proof in [2] is rather complicated 
and uses some special topological methods (rational models ets.).
We propose much more simpler (and purely algebraic) method to prove the main result from [2]. 
Moreover we prove that for Example 2 from [1] there is no compact form too. But it appears that some 
modification  of the construction of the Example 2 from [1] gives some other example of a solvable Lie 
group $R^{\prime}$ with a compact form  $M=R^{\prime}/\Gamma$.

The author express his thanks to A.Tralle for giving information about some modern problems
in symplectic topology.

Firstly we are going to describe the Lie algebras from the Examples 2 and 3 in [1].

Example 2. Here a Lie algebra $L(R)$ equals to $\text{Span}(A, B, X_1, X_2, X_3, Z_1, Z_2, Z_3)$,
where

$$[X_2, X_3]=2Z_1, [X_1, X_3]=Z_2, [X_1, X_2]=-Z_3,$$
$$[A, X_1]=-X_1, [A, X_2]=-2X_2, [A, X_3]=3X_3,$$
$$[A, Z_1]=Z_1, [A, Z_2]=2Z_2, [A, Z_3]=-3Z_3. $$

We have $L(R) = A + U$, where 

$$A = \text{Span}(A,B), U = \text{Span}(X_1,X_2,X_3,Z_1,Z_2,Z_3).$$

Example 3. Here $L(R)=\text{Span}(A,B,X_1,Y_1,Z_1,X_2,Y_2,Z_2)$, where

$$[X_1, Y_1]=Z_1, [X_2, Y_2]=Z_2,$$
$$[A, X_1]=X_1, [A, X_2]=-X_2,$$
$$[A, Y_1]=-2Y_1, [A, Y_2]=2Y_2,$$
$$[A, Z_1]=-Z_1, [A, Z_2]=Z_2 $$

Here we have $L(R)=A + U$, where 

$$A = \text{Span}(A,B), U=\text{Span}(X_1, Y_1, X_2, Y_2, Z_1, Z_2).$$

   Now we give some general considerations about Examples 2 and 3 from [1].

In both Examples the Lie group $R$ (simply connected) has a form $\bold R \times (\bold R \times _\phi N)$,
where $R^\prime = \bold R \times _\phi N$ is a semi-direct product, corresponding to a homomorphism $\phi : 
\bold R \to \text{Aut} (N)$. In both this Examples $N$ are some six-dimensional 2-nilpotent Lie groups, $L(N)=U$. 
For the Example 2 the Lie algebra for Lie group $N$ is the free 2-nilpotent Lie algebra with three generators, 
for the Example 3 we have $N=N_3(\bold R)\times N_3(\bold R)$, where $N_3(\bold R)$ is the unique three-dimensional 
non-Abelian simply connected nilpotent Lie group ($N_3(\bold R)$ is isomorphic to the group of all
real unipotent 3-matrices). In both cases the homomorphism $\phi$ has only real characteristic
numbers and therefore Lie groups $R$ are triangular (or, in other terms, completely solvable);
about these groups see, for example,[3]. For the triangular Lie groups there is some very important 
property: the algebra of cohomologies $H^*(L(R),\bold R)$  for the corresponding Lie algebra $L(R)$ is 
isomorphic to $H^*(R/\Gamma,\bold R)$ (where $\Gamma$ is a lattice in R, i.e. the discrete subgroup in 
$R$ with a compact factor-space $R/\Gamma$) [5].

  Also in both Examples we have $[R,R]=N$. We need some general result about the lattices in
triangular Lie groups. This result can be found in [6], but we prefer to give here some explicit 
(and more easier than in [6])  proof.

\proclaim{Proposition 1} Let $R$ be a triangular Lie group and $\Gamma$ be a lattice in $R$.
Then the commutator $[\Gamma,\Gamma]$ of $\Gamma$ is a lattice in $[R,R]$. In particular,
$\Gamma \cap [R,R]$ is a lattice in $R$.
\endproclaim

\demo{Proof} We may suppose that $R$ is simply connected. Consider $N=[R,R]$, $N$ is a 
simply connected nilpotent Lie group, therefore we have a natural structure of an algebraic 
group on $N$ (see,for example [7,4]). Let us consider the algebraic closure $N_1$ of a  subgroup 
$[\Gamma, \Gamma]$ in $N$, this subgroup is a lattice in $N_1$. It is clear that $\Gamma$ 
normalizes $N_1$ (because $\Gamma$ normalizes $[\Gamma,\Gamma]$). We have the following simple 
Lemma (see [8]):

\proclaim{Lemma 1 [8]} Let $R$ be a triangular connected Lie group and $F$  be a connected Lie 
subgroup in $R$. Then the normalizer $N_R (F)$ of $F$ in $R$ is connected.
\endproclaim

By the way, this Lemma shows that the condition 1). in the main result of [9](Theorem 1 - proof of 
the Benson-Gordon conjecture  in some very special case) may be excluded. But the second condition in 
this theorem are extremely strong and therefore there are a only few situations where this theorem may 
be applied.

   Due to this Lemma 1 we see that the normalizer $N_R(N_1)$ of $N_1$ in $R$ is connected ($N_1$ is 
connected as an algebraic subgroup in the connected nilpotent algebraic group $N$). But $N_R(N_1)$ 
contains $\Gamma$ and it is easy to understand that $N_R(N_1)$ is necessarily equals to $R$.

  Now we consider $R_1=R/N_1$,$\Gamma_1=\Gamma /[\Gamma,\Gamma]$,where $\Gamma_1$ is an Abelian lattice in
the triangular Lie group $R_1$. It is easy to understand (due to triangularity of $R_1$) that $R_1$ must be
Abelian too. Consequently $N_1 \supset [R,R]$, therefore $N_1=[R,R]$. We get that $[\Gamma,\Gamma]$ is a 
lattice in $R$. The intersection $\Gamma \cap [R,R]$ is a discrete subgroup in $[R,R]$ and contains
$[\Gamma,\Gamma]$, therefore $\Gamma \cap [R,R]$ is the lattice in $[R,R]$ too.
\enddemo

We continue our consideration of Examples 2 and 3 from [1]. For the nilpotent Lie group $N$, mentioned above, 
we have due to Proposition 1 that $\Gamma \cap N$ is a lattice in $N$ (in both Examples).

  It is well known (see, for example [7]) that if a nilpotent Lie group $N$ has a lattice, then the corresponding 
Lie algebra has a rational structure (i.e. it's structure constants in an appropriate basis are rational). If we
consider a natural structure of an algebraic group on $N$ (we suppose that $N$ is simply connected) then in 
this case $N$ will be defined over $\bold Q$ and the lattice in $N$ is commensurable with a subgroup $N_{\bold Z}$
of integer points of $N$.

 In [10] all nilpotent Lie algebras of dimension no more than 6 are classified (over an arbitrary field of 
characteristic 0). In particular we get a classification of the rational Lie algebras up to dimension 6.
This classification may be considered as a classification of the lattices (up to commensurability) in
real nilpotent simply connected Lie groups $N$. In is interesting to mention that in some such $N$ there are
an infinite series of non commensurable lattices (it is true, for example, for $N_3 \times N_3$).

 Let us consider now the Example 3 from [1] in details.  Here $N=N_3 \times N_3$. We have $L(R)=\bold R \oplus
L^\prime$, where $\bold R= \text{Span} (B), L^\prime=\bold R +_{\phi} U$ (here $\bold R= \text{Span}(A),U=L(N)$. 
An  action of $A$ on $U$ is defined by the matrix $\Phi = \text{diag}(1, -2, -1, -1, 2, 1)$. Let us suppose 
that there is some lattice $\Gamma$ in $R$.Then by virtue of Proposition 1 there is a lattice in $N$. 
Therefore on $N$ (and on the corresponding Lie algebra $L(N)$) it must be some rational structure. It is easy 
to understand that there is a lattice $\Gamma^\prime =\bold Z +_\psi (\Gamma \cap N)$ in $R^\prime$. Let $\gamma$ 
be a generator of $\bold Z$ in this decomposition. Then the  action of $\gamma$ on $L(N)$ is equal to the action 
of $C \cdot J$, where $C=\text{exp}(\Phi)t_0$ for some $t_0 \in \bold R$ and $J$ is some unipotent matrix (in fact 
$J$ is one of the elements of the adjoint Lie group for $N$).  The characteristic numbers of $C$ are $z,z,1/z,1/z,z^2,1/z^2$, where $z=exp(t_0), t_0\ne 0$. 

An action of $C$ on $L(N)$ induces an action on $L_{ab}=L(N)/[L(N),L(N)]$ - abelianization of $L=L(N)$.
For the lattice $D=\Gamma \cap N$ its intersection $D \cap [N,N]$ is a lattice in $[N,N]$ (see [7]),
consequently $D/D \cap [N,N]$ is a lattice in Abelian Lie group $N/[N,N]=\bold R^4$. The action of $C$ 
on $[N,N]$ preserves a lattice (which is isomorphic to $\bold Z^2$. The characteristic numbers of this 
action are $z,1/z$, therefore $z+1/z=n$ (as the trace of matrix) for some $n \in \bold N$.

Let us consider the rational structure (i.e. structure of rational Lie algebra)on $L(N)$, which is 
corresponding to the lattice $\Gamma \cap N$.

\proclaim{Lemma 2} Let $L$ be a Lie algebra over $\bold Q$, for which $L \otimes \bold R = 
n_3(\bold R)\times n_3(\bold R)$. Then $L$ is isomorphic (over $\bold Q$) to some Lie algebra
$L(p_1,p_2,  \dots p_k)$, where $p_i$ are the some pairwise different primes or zero, which defined 
by such commutator relations: 

$$[X_1,X_2]=X_5, [X_2, X_4]=X_5, [X_1,X_4]=X_6, [X_2, X_3]= q \cdot X_6$$, where $q=p_1 \cdot p_2 \cdot \dots 
\cdot p_k$.
\endproclaim

\demo{Proof} We use the classification of six-dimensional Lie nilpotent algebras over field $K=\bold Q$ 
(see [10]). Our lattice $D$ is 2-nilpotent and the rangs of its center and commutator equal to 2 (is follows 
from the corresponding properties of $N$ - see above). There  are only a few Lie algebras with such properties in 
the list from [10]. They are (we use the notation from [4]):

$g_{3,1}\oplus g_{3,1}$ (where $g_{3,1}=n_3$), $g_{5,1}\oplus \bold Q, g_{6,4}, g_{6,5}$. It is easy 
to verify that for these algebras tensored by $\bold R$ only $g_{3,1}\oplus g_{3,1}$ and $g_{6,5}$ are isomorphic 
to $n_3(\bold R)\times n_3(\bold R)$. For $g_{6,5}$ the commutator relations are:

$$[X_1,X_2]=X_5, [X_2, X_4]=X_5, [X_1,X_4]=X_6, [X_2, X_3]= q \cdot X_6$$,

where $q \ge 0$ belongs to $\bold Q / \bold Q^2$ (for $q=0$ we get $g_{3,1}\oplus g_{3,1}$). Therefore
we may suppose that  $q=p_1 \cdot p_2 \cdot \dots \cdot p_k$.
\enddemo

For $q=2$ some construction of the corresponding rational Lie algebra can be found in [10].

  Now we consider Lie algebras $g_{6,5}$ for $q >0$. It is easy to calculate the group of their automorphisms.
There are two special one-dimensional groups of these automorphisms - $A_\alpha, B_\beta$:

$$ A_\alpha : X_1 \to X_1,X_2 \to X_2,X_3 \to \alpha X_3, X_4 \to \alpha X_4, X_5 \to \alpha X_5, 
X_6 \to \alpha X_6,$$

$$ B_\beta : X_1 \to \beta X_1,X_2 \to \beta X_2,X_3 \to X_3, X_4 \to X_4, X_5 \to \beta X_5, 
X_6 \to \beta X_6,$$

where $\alpha, \beta \in \bold R \setminus \{ 0 \}$.

 It is easy to proof that $\text{Aut}(g_{6,5})=(A \cdot B)\cdot U$, here $F=A \cdot B$ is the reductive part
(an Abelian two-dimensional group) of $\text{Aut}(g_{6,5})$ and $U$ -  its unipotent radical ($\text{dim} U =7$).
A matrix realization of $F$ can be written in a such form:

$$F= \{ \text{diag}(\beta,\beta, \alpha,\alpha,\alpha \beta, \alpha \beta),$$ 

the latter two elements $\alpha \beta$ are corresponding to the action on the center of our Lie algebra.

  We have a rational structure on $L(N)$ due to lattice in $N$ (see above). In virtue of Lemma 2 this 
structure is isomorphic to $g_{6,5}$ for some $q \ge 0$. Our first case will be $q >0$. The lattice in $R^\prime$ 
gives as an automorphism of $L(N)$ with the characteristic numbers  $z,z,1/z,1/z,z^2,1/z^2$, where $z=exp(t_0)$. 
We have $z>0$. The only case when the set of these numbers  $z,z,1/z,1/z,z^2,1/z^2$ equals to a set 
$(\beta,\beta, \alpha,\alpha,\alpha \beta, \alpha \beta)$ is when $z=1$. But $z=exp(t_0)$ and $t_0 \ne 0$.
We have a contradiction.

 Now we consider the case when $q=0$, here $L=L(N)=n_3 (\bold Q) \oplus n_3 (\bold Q)$. In this case 
the group of automorphisms has a  decomposition $GL_2(\bold Q)\times GL_2(\bold Q) \cdot U$, here $U$ is the 
nilradical ($\text{dim}U=8)$, two $GL_2(\bold Q)$ are corresponding to the groups of linear transformations of $n_3/[n_3,n_3]$. Let us consider $L_{ab}=L/[L,L]=\bold Q^4$ and the induced action of $C$ on this vector space. 
The characteristic numbers of this action are $z,1/z,z^2,1/z^2$. We know that $z$ is real and positive. 
If $z \ne 1$ then all four numbers mentioned above are distinct. Therefore is this case there are only three 
decompositions of $L_{ab}$ into a direct sum of two two-dimensional subspaces. If $e_1,e_2,e_3,e_4$ are proper 
vectors, corresponding to $z,1/z, z^2, 1/z^2$, then these decompositions are $\text{Span}(e_1,e_2)\oplus \text{Span}(e_3,e_4), \text{Span}(e_1,e_3) \oplus \text{Span}(e_2,e_4), \text{Span}(e_1,e_4) \oplus 
\text{Span}(e_2,e_3)$. We are going to show that the first and the last decompositions over $\bold Q$ are impossible.

 Let us suppose that the first decomposition is defined over $\bold Q$ (for the rational structure on $L=L(N)$).
Then $z+z^2=k,1/z+1/x^2=l$ for some natural $k,l\in \bold N$. From this relations it follows that $z= (kl-1)/l$,
i.e. $z$ must be rational. But we have also $z+1/z \in \bold N$.It is easy to understand that it is possible
only if $z=1$. We know that $z \ne 1$, therefore this decomposition into the direct sum is impossible. Analogously 
we can prove that the third decomposition is impossible too. 

  As we consider the case $q=0$, when the rational structure on $L(N)$ is isomorphic to $n_3 \oplus n_3$, it is 
necessary to  have some decomposition for $n_{ab}$ (which is corresponding to the decomposition of $L$ into a 
direct sum due to its definition). Moreover, this decomposition of $L_{ab}$ must generate some  decomposition of 
$L$ into a direct sum. As it is proved above,we have only one variant of decomposition:$L_{ab}=\text{Span}(e_1,e_3) 
\oplus \text{Span}(e_2,e_4)$. But $[e_1,e_2]=[e_2,e_4]=0$ (due to commutator relations), therefore this decomposition can not correspond to the decomposition of $L$ into a direct sum of ideals ($e_1,e_3$ cannot generate such ideal, also for $e_2,e_4$). Once again we get a contradiction. 

In all subcases we have the contradictions. Therefore there is no lattices in the Lie group $R$ from the Example 
3 in [2]. 

  Now we are proceed to the Example 2 from [1]. Here we are going to prove  the nonexistence of any lattice in the corresponding Lie group $R$ too. But for this example we'll give below some modification $R^\star$ of $R$ with the 
same geometrical and cohomological properties as $R$ and with some lattice.

   The proof of nonexistence of a lattice for $R$ from the Example 2 is analogous but much more shorter that in the case 
of the Example 3. Here we have corresponding Lie group $R^\prime = \bold R \cdot_\phi N$. We suppose that there is    
a lattice in $R$, then we have a lattice $\Gamma = \bold Z \cdot D$ in $R^\prime$, where $D=\Gamma \cap N$ is a 
lattice in $N$. The Lie algebra of $N$ is the free 2-nilpotent Lie algebra with 3 generators, it is $g_{6,3}$  
over $\bold Q$ in [10]).

   We have $L=V+\Lambda^2V, \text{dim}V=3, Z(L)=\Lambda^2V$. The action on $L$ of the generator $B \in \bold R$
in the decomposition of $L(R)$ has the characteristic numbers $-1,-2,3,1,2,-3$ (it follows from the
definition of $L$). Therefore the action of the generator $\gamma$ of $\bold Z$ in the decomposition of $\Gamma$ has the characteristic numbers $1/z,1/z^2,z^3, z,z^2,1/z^3$, where $z=exp(Bt_0),t_0 \ne 0$. As the commutator is invariant 
under this action of $\gamma$, for the corresponding rational structure on $L$ we have two actions, defined over 
$\bold Q$ - on $V$ and $Z(n)=\Lambda^2$ (three-dimensional vector spaces).  The corresponding characteristic numbers 
for the action on $V$ are $1/z,1/z^2,z^3$, therefore, in particular, for their sum (=trace) we have $1/z+1/z^2+z^3= 
m$ for some $m \in N$.  analogously for the action on $\Lambda^2V$ the characteristic numbers are $z, z^2,1/z^3$, hence 
$z+z^2+1/z^3=n$ for some $n \in N$. We get two equations, they may be rewritten in a form of a polynomial system of equations:

$x^5-mx^2+x+1=0,$ 

$x^5+x^4-nx+3+1=0.$

   We are going to solve this system. For this we use the method of finding the Groebner Basis for the systems of polynomial equations [11]. With aid of the computer program Maple V we get such set of polynomials (they generate the 
same polynomial ideal as the equations from our system) as the Groebner Basis for our system:

$-7m^2 + m^3 - m^4 - m^5 - 13mn - m^2n + 5m^3n - 3m^4n - 7n^2 - mn^2 + 10m^2n^2 + n^3 + 5mn^3 + m^3n^3 - n^4 - 3mn^4 
- n^5 ,$ 

$-28m + 4m^2-4m^3-4m^4 - 20n - 16mn + 29m^2n - 8m^3n - m^4n - 8n^2 + 8mn^2 + 17m^2n^2-m^3n^2-5mn^3+12m^2n^3+2m^3n^3+ m^4n^3-10n^4- 12mn^4 - -2m^2n^4+3m^3n^4-6n^5-12mn^5+ m^2n^5-6n^6-m^2n^6+2n^7+mn^7+12nx+8n^2x+ 23n^3x + 9n^4x + 12n^5x+n^7x-n^8x,$ 

$-52m+6m^2+4m^3-6m^4-20n-68mn+63m^2n - 14m^3n-m^4n-52n^2+26mn^2+24m^2n^2+3m^3n^2+3m^4n^2+16n^3- 13mn^3 - 4m^2n^3 + 9m^3n^3-
- 18n^4-39mn^4+4m^2n^4-20n^5-mn^5-3m^2n^5 +6n^6+ 3mn^6+40mx+68nx+18n^2x+33n^3x+ 22n^4x+2n^5x+4n^6x-3n^7x,$  

$-20-26m-7m^2-8m^3-3m^4-20n+36mn+9m^2n-7m^3n+2m^4n+24n^2+33mn^2-18m^2n^2+4m^3n^2- m^4n^2 + 18n^3- 14mn^3+ 8m^2n^3 - 3m^3n^3--14n^4+8mn^4-3m^2n^4+10n^5+2mn^5+m^2n^5-2n^6- mn^6-20x+4nx-31n^2x+9n^3x-14n^4x+6n^5x-3n^6x+n^7x-20x^2-10nx^2-10n^2x^2,$

$40 - 52m+6m^2+4m^3-6m^4-20n-68mn+63m^2n-14m^3n-m^4n-52n^2+26mn^2+24m^2n^2 + 3m^3n^2 + 3m^4n^2 + 16n^3 - 13mn^3 - 4m^2n^3+
+ 9m^3n^3 - 18n^4 - 39mn^4 + 4m^2n^4 - 20n^5 - mn^5 - 3m^2n^5 + 6n^6 + 3mn^6 + 68nx + 18n^2x + 33n^3x + 22n^4x + 2n^5x + 4n^6x - 3n^7x + 40nx^2-40x^3.$

 Let us consider the corresponding polynomial equations.   First equation is the resultant of our system. Next two equations are linear in $x$. If in one of these two equations the coefficient for $x$ is nonzero, we get that $x$ is rational. But a rational root of our initial equation$x^5-mx^2+x+1=0$ have to divide 1, therefore this root equals to 1 (we remind that $x > 0$).  For the second equation the coefficient for x equals to $12n+8n^2+23n^3+9n^4+12n^5+n^7-n^8$, the only root in $\bold N$ of this polynomial is $n=3$. For the coefficient from the second equation we get  $-120+40m=0$ (using $n=3$), therefore m=3. Our initial system for $m=n=3$ has only one solution $x=1$. Therefore we get $x=1$, but it is a contradiction and there is no lattices in $R$ from Example 2 in [1].

  Let us consider some variant of the Example 2 from [1]. Let $L$ be a Lie algebra with basis 

$$A, B,X_1,X_2, X_3, Z_1,Z_2, Z_3$$ 

and commutator relations

  $$[X_1, X_2]= -Z_3, [X_1, X_3]=Z_2,[X_2, X_3]=2X_1$$
  $$[A,X_1]=\lambda_1 X_1,[A,X_2]=\lambda_2X_2,[A,X_3]=\lambda_3 X_3,$$
  $$[A,Z_1]=
(\lambda_2+\lambda_3)Z_1,[A,Z_2]=(\lambda_1+\lambda_3)Z_2, [A,Z_3]=(\lambda_1+\lambda_2)Z_3.$$ 

for some $\lambda_i \in \bold R$. Following [1] we set

$$\omega=\alpha \wedge \beta +\mu_1 \wedge \zeta_1+\mu_2 \wedge \zeta_2+\mu_3 \wedge \zeta_3$$

where $\alpha,\beta,\mu_1,\mu_2,\mu_3,\zeta_1,\zeta_2,\zeta_3$ - the dual basis in $L$. 

  It is easy to check that all the properties of $L$ and $ \omega$ from [1] are true for our more general Lie algebra $L$ if $\lambda_i \ne 0, i=1,2,3$ and $\lambda_1+\lambda_2+\lambda_3=0$. Therefore as in [1] we get the Hard Lefschetz property for our $L$. Also the cohomology ring for $L$ is isomorphic to the cohomology ring for the K\"{a}hler manifold $T^2 \times \bold{CP}^3$. We are going to prove that for some triples $\lambda_1,\lambda_2,\lambda_2$ we can construct a lattice in a simply connected triangular Lie group $R$ corresponding to our Lie algebra $L$.

 Let $x^3-px^2+qx-1$ be a polynomial with integer coefficients and three real roots. For example, the equation 
$x^3-5x^2+6x-1=0$ has the roots $x_1\approx 0.198, x_2 \approx 1.555, x_3 \approx 3,247$.
We take $\lambda_i=\text{ln}(x_i)$,in our example $\lambda_1 \approx -1.619,\lambda_2 \approx 0.441,\lambda_3 \approx 1.777$.

For these $\lambda_i$ the action of the matrix $C=\text{diag}(exp(\lambda_1, \lambda_2,\lambda_3)$ on $V$ is conjugate to some action from $\text{GL}(3,\bold Z)$. The corresponding action on $\Lambda^2V$  is integer-valued too. Therefore there is a lattice $D$ in $N$ which is invariant under the action of some element $\gamma$ from the group of automorphisms. We set   $\Gamma_1 = \bold Z \times _{\phi} D$, it is a lattice in $R_1$. The group $\Gamma=\Gamma_1 \times \bold Z$ is a lattice in $R=R_1 \times \bold R$. Therefore we get a lattice in our Lie group $R$ (which is slightly different from $R$ 
in Example 2 [1]). The compact symplectic manifold $M=R/\Gamma$  gives us an example which failed to be found in [1]. For this $M$ all the cohomological properties are the same as for K\"{a}hler manifold. Moreover,the (nilpotent) minimal model for $M$ is the same that for the K\"{a}hler manifold $T^2 \times \bold{CP}^3$ (as in [1]). It is not known if 
there is a complex structure on this $M$. Is it is really no complex structure on this $M$, therefore there is no K\"{a}hler structure on it. In this case we'll find that there is no cohomological invariants which can characterize 
K\"{a}hler solvmanifolds. Some other examples of such kind are, as indicated in [1], six-dimensional manifolds $Sol 
\times T^3$ and $Sol \times Sol$,where $Sol$ is constructed in [1],Example 1,

\enddocument